\documentclass[12pt,reqno]{amsart}

\newtheorem{theorem}{Theorem}[section]

\newtheorem{lemma}[theorem]{Lemma}

\theoremstyle{definition}
\newtheorem{definition}[theorem]{Definition}

\newcommand{\up}[1]{{\ensuremath{^{<{#1}>}}}}

\newcommand{\dnm}[1]{\ensuremath{_{<{#1}>}}}

\newcommand{\doum}[1]{\ensuremath{_{<\! \! \!<{#1}>\! \! \!>}}}
\newcommand{\main}{\ensuremath{ \mathrm{K}_i(V)}}
\newcommand{\sub}{\ensuremath{ \mathrm{D}_i(V)}}

\newcommand{\subj}{\ensuremath{ \mathrm{D}_j(V)}}

\begin{document}
\title[A combinatorial proof of persistence theorem] {A combinatorial proof of Gotzmann's persistence theorem for monomial ideals}
\author{Satoshi Murai}
\maketitle
\vspace{-1pt}
\begin{abstract}
Gotzmann proved the persistence for minimal growth for ideals.
His theorem is called Gotzmann's persistence theorem.
In this paper,
based on the combinatorics on binomial coefficients,
a simple combinatorial proof of Gotzmann's persistence theorem in the special case of monomial ideals is given.
\end{abstract}
\vspace{-1pt}
\section*{Introduction}
Let $K$ be an arbitrary field, $R=K[x_1,x_2,\dots,x_n]$  the polynomial ring with deg$(x_i)=1$ for $i=1,2,\dots,n$.
Let $M$ denote the set of variables $\{ x_1,x_2,\dots,x_n\}$, 
$M^d$ the set of all monomials  of degree $d$, where $M^0=\{1\}$
, and  $\overline{M_i}=M\setminus \{x_i\}$.
For a monomial $u\in R$ and for a subset $V\subset M^d$,
we define
$uV=\{uv|v\in V\}$ and
$MV=\{x_iv |v\in V,\ i=1,2,\dots,n\}$.
For a finite set $V\subset M^d$,
we write $|V|$ for the number of the elements of $V$. 
Let $\gcd(V)$ denote the greatest common divisor of the monomials belonging to $V$.

Let $n$ and $h$ be positive integers.
Then $h$ can be written uniquely in the form,
called the $n$\textit{th binomial representation of} $h$,
$$h={{h(n) +n}\choose n}+{{h(n-1)+n-1}\choose n-1}+\dots+{{h(i)+i}\choose i} ,$$
where $h(n)\geq h(n-1)\geq \dots\geq h(i)\geq 0,\ i\geq 1$.
See [3, Lemma 4.2.6].

Let ${h(1)\choose s(1)}+{h(2)\choose s(2)}+\dots+{h(i)\choose s(i)} $ be a sum of binomials,
where\ $h(j)\geq s(j)$ for any $j=1,2,\dots,i  $.
Then we define\begin{eqnarray*}
\bigg\{{h(1) \choose s(1)}+\dots+{h(i)\choose s(i)} \bigg\}^{[+1]}&=&{h(1)+1 \choose s(1)}+\dots+{h(i)+1\choose s(i)}.\\ 
\end{eqnarray*}
Let $h={h(n) +n\choose n}+\dots+{h(i)+i\choose i} $ be the $n$th binomial representation of $h$.
We define \begin{eqnarray*}
h^{ <  n  >}&=&{h(n)+n+1 \choose n}+\dots+{h(i)+i+ 1\choose i} ,\\
h_{<  n >}&=&{h(n)+n \choose n - 1}+\dots+{h(i)+i\choose i - 1},\\
h_{<\!\!\!<  n>\!\!\!>}&=&{h(n)+n - 1 \choose n - 1}+\dots+{h(i)+i - 1\choose i- 1} ,\\
\end{eqnarray*}
and  set $0^{< n>}=0_{<  n>}=0_{<\!\!\!<  n>\!\!\!>}=0$, $1^{<  0> }=1_{<\!\!\!<   0>\!\!\!>}=1$ together with $1_{< 0>}=0$.

The  inequality (\ref{mb}) below was proved by F. H. S. Macaulay.
See also \cite{CM} and \cite{IS} for further infomation.
Let $V$ be a set of monomials of same degree.
Then one has 
\begin{eqnarray}
|MV|\geq |V|\up{n-1}.\label{mb} 
\end{eqnarray}
In 1978, Gotzmann \cite{G} proved so-called persistence theorem.
In the special case of monomial ideals, the persistence theorem says that
\begin{theorem}[\textbf{Persistence Theorem for monomial ideals}]\label{main}
Let $V$ be a set of monomials of  degree $d$.
If $|MV|=|V|\up{n-1}$,
then $|M^{i+1}V|=|M^iV|\up{n-1}$ for all $i\geq 0$.
\end{theorem}

Let $A=(a_1,a_2,\dots,a_n)$ and $B=(b_1,b_2,\dots,b_n)$ be elements of $\mathbb{Z}_{\geq 0}^n$.
The \textit{lexicographic order} on  $\mathbb{Z}^n$ is defined by $A<B$ if the leftmost nonzero entry of $B-A$ is positive.
Moreover, the lexicographic order on monomials of the same degree is defined by
${x_1}^{a_1}{x_2}^{a_2}\dots{x_n}^{a_n}<{x_1}^{b_1}{x_2}^{b_2}\dots{x_n}^{b_n}$ if $A<B$ on $\mathbb{Z}_{\geq 0}^n$.

Let $V$ be a set of monomials of degree $d$.
\begin{itemize}
\item[(i)]
$V$ is called a \textit{Gotzmann set} if  $V$ satisfies $|MV| =|V|\up{n-1}$.
\item[(ii)]
$V$ is called a \textit{lexsegment set} if $V$ is a set of first $|V|$ monomials in lexicographic order.
Denote the lexsegment set $V$ of $K[x_1,\dots,x_n]$ in degree $d$ with $|V|=a$ by $Lex(n,d,a)$.
\end{itemize}

It is known that lexsegment sets are  Gotzmann sets.
See \cite[\S 4.2]{CM} or \cite{IS}.
Also, in \cite{S} we determined all integers $a>0$ such that every Gotzmann set with $|V|=a$ and with $\gcd(V)=1$ is lexsegment up to permutation of variables.  
Related works of Gotzmann's theorem   were done by A. Aramova, J. Herzog and T. Hibi  \cite{AH}.
They proved Gotzmann's theorem for exterior algebra.
In addition, 
Z. Furedi and J. R. Griggs \cite{FG} determine all integers $a>0$  such that every squarefree Gotzmann set with $|V|=a$ is squarefree lexsegment up to permutation. 

The inequality (\ref{mb})  and Theorem \ref{main} are true for more general case.
They need not to be restricted to monomial case.
Gotzmann  \cite{G} proved the persistence for minimal
growth of the Hilbert function of a homogeneous ideal
(see \cite[Theorem C.17]{IS}).
M. Green refined Gotzmann's proof
 (see \cite[Theorem 4.3.3]{CM}).
Green also give a simple proof in \cite[Theorem 3.8]{Gr} using  generic initial ideals.
On the other hand, in the special case of monomial ideals,
in \cite{G} Gotzmann  proved the persistence theorem easier than general case 
using his version of the theory of Castelnuvo--Munford regularity.
All of these proofs are completely algebraic.
In the present paper we will give a combinatorial proof of  persistence theorem for monomial ideals.
The advantage of our proof is that we only use the combinatorics on binomials.

In \S 1, we will prepare some lemmas about binomial representations.
In \S 2, we will give a combinatorial proof of persistence for monomial ideals. 

\section {Binomial representations}
In this section we consider some properties about binomial representation and combinatorics
which will be used in the main proof.
\begin{definition}\label{amari}
Let $h$ be a positive integer and $h=\sum_{j=i}^n {h(j) +j \choose j}$ the $n$th binomial representation of $h$.
Let $\alpha =$max$\{0, $max$\{\alpha\in \mathbb{Z}|h-{{\alpha + n}\choose n}>0\}\}$.
We denote $h-{{\alpha +n}\choose n}$ by $\bar{h}^{(n)}$, in other words,
\begin{itemize}
\item[(i)] if $h=1$, then $\bar{h}^{(n)}=0$;
\item[(ii)] if $h>1$ and $i=n$, then $\bar{h}^{(n)}={{h(n)+n-1} \choose{n-1}}$;
\item[(iii)] if $h>1$ and $i<n$, then $\bar{h}^{(n)}=\sum_{j=i}^{n-1} {{h(j)+j}\choose j}$.
\end{itemize}
This constraction says  $\bar{h}^{(n)}\leq {{\alpha +n}\choose n-1}$ and
 $h\up{n}={{\alpha+n}\choose n}\up{n} +{\bar{h}^{(n)}}\up{n-1}$.
Furthermore, if $h>1$ then   $\bar{h}^{(n)}\geq1$.

Firstly, we introduce some easy and fundamental properties.
\end{definition}

\begin{lemma}[{[3 Lemma 4.2.7]}] \label{b-unique}
Let $a=\sum_{k=i}^n {{h(k)}\choose k}$ and $a'=\sum_{k=j}^n {{h'(k)}\choose k}$ be the binomial representations.
Then one has $a<a'$
if and only if
\[(h(n),h(n-1),\dots,h(i),0,\dots,0)<(h'(n),h'(n-1),\dots,h'(j),0,\dots,0)\]
in the lexicographic order on $\mathbb{Z}^n$ .
\end{lemma}

\begin{lemma}\label{apply}
Let $h$ and $n$ be  integers with $h\geq 0$ and $n>0$.
Then,  for any integer $1\leq \alpha \leq h$, one has
\[{{h+n}\choose n}={{\alpha-1+n}\choose n}+{{\alpha+n-1}\choose n-1}+{{\alpha+1+n-1}\choose n-1}+\dots+{{h+n-1}\choose n-1}\]
and
\[{{h+n}\choose n}^{[+1]}=\bigg\{{{\alpha-1+n}\choose n}+{{\alpha+n-1}\choose n-1}+{{\alpha+1+n-1}\choose n-1}+\dots+{{h+n-1}\choose n-1}\bigg\}^{[+1]} .\]
\end{lemma}
\begin{proof}
Use ${{h+n}\choose n }={{h-1+n}\choose n}+{{h-1+n}\choose n-1}$ to the leftmost binomial coefficient repeatedly,
then  we have
\begin{eqnarray*}
{{h+n}\choose n }&=&{{h-2+n}\choose n }+{{h-1+n-1}\choose n-1 }+{{h+n-1}\choose n-1 }\\
&\vdots& \\
&=&{{\alpha-1+n}\choose n}+{{\alpha+n-1}\choose n-1}+{{\alpha+1+n-1}\choose n-1}+\dots+{{h+n-1}\choose n-1},
\end {eqnarray*}
as desired.
\end{proof}

\begin{lemma}\label{change}
Let $h$ and $n$ be positive integers. Then,
$$h\up{n}=h+h\dnm{n}.$$
\end{lemma}
\begin{proof}
Let $h=\sum_{j=i}^n {{h(j)+j}\choose j}$ be the $n$th binomial representation of $h$.
Since ${{h+n}\choose n }={{h-1+n}\choose n}+{{h-1+n}\choose n-1}$, one has
$$h+h\dnm{n}=\sum_{j=i}^n \{{{{h(j)+j} \choose j}+{{h(j)+j}\choose j-1}}\}
=\sum_{j=i}^n {{h(j)+j+1} \choose j}=h\up{n},$$
as desired.\end{proof}

Next, we introduce some lemmas which will be used in the main proof.

\begin{lemma}\label{plus}
Let $a$, $b$ and $n$ be positive integers. One has
$$a\up{n}+b\up{n}>(a+b)\up{n}.$$
\end{lemma}
\begin{proof}
Assume $n\geq 2$.
Then we can take  $d$ with $|M^d|>a+b$.
Let $V_a=Lex(n,d,a)$, $V_b=Lex(n,d,b)$ and $u$  the minimal element of $V_a$ in the lexicographic order.
Let $V=x_1^{d+1}V_a\cup ux_n V_b$.
Since $ux_1^{d+1}>ux_1^dx_n$,
$x_1^{d+1}V_a \cup ux_n V_b$ is disjoint union if $n\geq2$.
Since $x_1^{d+1}x_nu\in Mx_1^{d+1}V_a\cap Mux_nV_b$,
we have $Mx_1^{d+1}V_a\cap Mux_nV_b\ne \emptyset$.
By (\ref{mb})  for any positive integer $n\geq2$,
we have
$$(a+b)\up{n-1}\leq |MV|<|MV_a|+|MV_b|=a\up{n-1}+b\up{n-1},$$
as desired.
\end{proof}

\begin{lemma}\label{combinatrix}
Let $a,b,c$ and $\alpha$ be positive integers.
If ${{\alpha+n} \choose n} +a=b+c$ and \linebreak
$a,b,c<{{\alpha+n}\choose n}$,
then one has  $${{\alpha+n} \choose n}\up{n} +a\up{n}\leq b\up{n}+c\up{n}.$$
Especially, if  ${{\alpha+n} \choose n}\up{n} +a\up{n}=b\up{n}+c\up{n}$, then we have
\begin{eqnarray}
\bigg\{{{\alpha+n} \choose n}\up{n}\bigg\}\up n +\{a\up{n}\}\up n 
=\{b\up{n}\}\up n+\{{c\up{n}} \}\up n. \label{second}
\end{eqnarray}
\end{lemma}
\begin{proof}
We use induction on $n$.
\medskip

\noindent \textbf{[Case I]} Let $n=1$.

In general, if $h$ is a positive integer, then $h\up{1}={{h+1}\choose 1}=h+1$.
Thus we have
${{\alpha+1} \choose 1}\up{1} +a\up{1} = b+1+c+1 =b\up{1}+c\up{1} $.
Thus we may assume $n>1$.

To prove Lemma \ref{combinatrix}, we claim the followings:
\medskip
\begin{itemize}
\item[\textbf{(\#\#)}]Let $h$ and $s$ be positive integers.
Assume Lemma \ref{combinatrix} is true in the case of $n=s$.
If ${{h+s}\choose s}=\sum_{i=1}^k h_i+c-d$,
 $0<h_i<{{h+s}\choose s} $ for $i=1,2,\dots,k$, $k\geq 1$ and ${{h+s}\choose s}>c>d\geq 0$, 
then one has
$${{h+s}\choose s}\up{s}\leq\sum_{i=1}^k h_i\up{s}+c\up{s}-d\up{s}.$$
Especially, if ${{h+s}\choose s}\up{s}=\sum_{i=1}^k h_i\up{s}+c\up{s}-d\up{s}$, then we have
$$\bigg\{{{h+s}\choose s}\up{s}\bigg\}\up s = \sum_{i=1}^k \{h_i\up{s}\}\up s+\{c\up{s}\}\up{s}-\{d\up{s}\}\up s.$$
\end{itemize}
\bigskip
We will prove the claim. 
Since Lemmas \ref {combinatrix} is true for $n=s$, we have \begin{eqnarray}
{{h+s}\choose s}\up{s}\leq \{\sum_{i=1}^k h_i \}\up{s} +c \up{s}-d\up{s}. \label{tuika}
\end{eqnarray}
Moreover, by Lemma \ref{plus}, we have
\begin{eqnarray} 
\{\sum_{i=1}^k h_i \}\up{s} +c\up{s} \leq \sum_{i=1}^k h_i\up{s} + c\up{s}. \label{tuiki}
\end{eqnarray}
Also, if  $k\geq 2$ then (\ref{tuiki}) is not equal. 
If $k=1$, then (\ref{tuika}) is of the form Lemma \ref{combinatrix}.
Thus by  an assumption we proved the claim (\#\#).
\medskip

We return to the proof of Lemma \ref{combinatrix}.
Let $a={{a(n)+n}\choose n}+\bar{a}^{(n)}$,
$b={{b(n)+n}\choose n}+\bar{b}^{(n)}$ and
$c={{c(n)+n}\choose n}+\bar{c}^{(n)}$ be the form of Definition \ref{amari}.
Let $\bar{a}=\bar{a}^{(n)}$, $\bar{b}=\bar{b}^{(n)}$ and $\bar{c}=\bar{c}^{(n)}$.
First, we note fundamental inequalities.
\begin{itemize}
\item[($\alpha$)] $a<b,a<c,\alpha>b(n),\ \alpha>c(n)\   \mathrm{ and }\  a(n)\leq c(n),$ 
\item[($\beta$)] $\bar{b}\geq1\  \mathrm{ and }\  \bar{c}\geq1,$
\item[($\gamma$)] $\bar{b}<{{\alpha+n-1}\choose n-1}\ \mathrm{ and }\ \bar{c}<{{\alpha+n-1}\choose n-1}.$
\end{itemize}
The inequality($\alpha$)  follows from the assumption.
We have the inequality($\beta$) since $1\leq a<b,c$.
By Definition \ref{amari}, we have
$\bar{b}\leq {{b(n)+n}\choose n-1} \leq {{\alpha+n-1}\choose n-1}$.
But if $\bar{b}= {{b(n)+n}\choose n-1}$, then $b(n)<\alpha -1$ since $b={{b(n)+1+n}\choose n}<{{\alpha +n}\choose n}$.
Thus we have the inequality($\gamma$).
Next, by Lemma \ref {apply}, we can write ${{\alpha+n}\choose n}$ and ${{c(n)+n}\choose n}$ as follows:\begin{eqnarray*}
{{\alpha+n}\choose n}= {{b(n)+n}\choose n}+\sum_{i=b(n)+1}^\alpha  {{i+n-1}\choose n-1};\\
{{c(n)+n}\choose n}= {{a(n)+n}\choose n} +\sum_{i=a(n)+1}^{c(n)}  {{i+n-1}\choose n-1}.
\end{eqnarray*}
Hence we substitute these equalities for  ${{\alpha+n} \choose n} +a=b+c$, then we have\begin{eqnarray}
\bigg \{ \sum_{i=b(n)+1}^\alpha {{i+n-1}\choose n-1} \bigg\}+\bar{a}=\bar{b}+\bar{c}+\bigg\{\sum_{i=a(n)+1}^{c(n)} {{i+n-1}\choose n-1} \bigg\}.\label{15}\end{eqnarray}
Furthermore,   ${{\alpha+n} \choose n}\up{n} +a\up{n}\leq b\up{n}+c\up{n}$  if and only if 
\vspace{-4pt}\begin{eqnarray}
&\bigg \{&\! \! \! \sum_{i=b(n)+1}^\alpha {{i+n-1}\choose n-1} \bigg\}^{[+1]}+\bar{a}\up{n-1} \nonumber \\
& &\leq \bar{b}\up{n-1}+\bar{c}\up{n-1}+\bigg\{\sum_{i=a(n)+1}^{c(n)} {{i+n-1}\choose n-1}\bigg\}^{[+1]}\label{16} .\end{eqnarray}
Instead of considering ${{\alpha+n} \choose n} +a=b+c$ and ${{\alpha+n} \choose n}\up{n} +a\up{n}\leq b\up{n}+c\up{n}$,
it is enough to consider (\ref{15}) and (\ref{16}).
We will consider two cases.
\medskip

\noindent \textbf{[Case II]} Let $\bar{c}\geq \bar{a}$ and $ n>1$.
We will prove that for $i=0,1,\dots,\alpha-(b(n)+1)$ ${{\alpha -i+n-1}\choose n-1}$ can be written
\begin{eqnarray}
{{\alpha-i+n-1} \choose n-1} =-d_i+\sum_{j=t_{i+1}}^{t_{i}-1} P_j+d_{i+1}, \label{X}
\end{eqnarray}
where $P_i={ i+n-1 \choose n-1}$ for $i = a(n)+1,\dots,c(n)$, $t_{i+1}<t_i\leq c(n)-i+2$,
 $0\leq d_i<P_{t_i-1}$
 together with $P_{c(n)+1}=\bar{c}$, $t_0=c(n)+2$, $d_0=\bar{a}$ and $d_{\alpha - b(n)}=\bar{b}$.

We use induction on $i$.
For $i=0$,
since $\bar{c}-\bar{a}<{{\alpha+n-1} \choose n-1}$,
there exists $t_1\leq c(n)+1 $ such that
$$\sum_{i=t_1}^{c(n)} {{i+n-1}\choose n-1}+\bar{c}-\bar{a}\leq {{\alpha+n-1} \choose n-1} < \sum_{i=t_1-1}^{c(n)} {{i+n-1}\choose n-1} + \bar{c}-\bar{a}.$$
Thus we have 
\begin{eqnarray*}
{{\alpha+n-1} \choose n-1} =\bar{c}-\bar{a}+\sum_{j=t_1}^{c(n)} P_j+d_1=-d_0+\sum_{j=t_1}^{c(n)+1} P_j+d_1\end{eqnarray*}
with $0\leq d_1<P_{t_1-1}$.
Assume we have the form (\ref{X}) for $i=0,\dots,s-1$.
By the assumption of induction and $\alpha >c(n)$ we have 
${{\alpha-s+n-1} \choose n-1}\geq {{c(n)-s+1+n-1} \choose n-1}\geq {{t_{s}-1+n-1} \choose n-1}=P_{t_{s}-1}$.
Thus ${\alpha -s +n-1\choose n-1} \geq -d_{s+1}+P_{t_s-1}$ and $t_{s+1}<t_s$.
By the same way of $i=0$,
we have (\ref{X}) for $i=s$.
Especially, if $s=\alpha -(b(n)+1)$, because of the equality (\ref{15}), we have
\begin{eqnarray}
{{b(n)+n}\choose n-1}=-d_s+\sum_{j=a(n)+1}^{t_{s}-1} {{j+n-1}\choose n-1}+\bar{b}.\label{19}\end{eqnarray}
Thus  each ${{\alpha -i+n-1}\choose n-1}$ have of the form (\ref{X}).

Equalities (\ref{X}) satisfies conditions of (\#\#).
By the assumption of induction of $n$, we have \begin{eqnarray}
 {{\alpha-i-1+n-1} \choose n-1}\up{n-1} \leq -d_{i+1}\up{n-1}+\sum_{j=t_{i+1}}^{t_{i}-1} P_j\up{n-1}+d_{i+2}\up{n-1}.\label{X'}
\end{eqnarray}
Summating  (\ref{X}) in both sides yields (\ref{15}), and
summating inequalities (\ref{X'}) in both sides yields (\ref{16}). 
Furthermore, (\ref{16}) is equal if and only if (\ref{X'}) are equal for all $i$.
Thus if (\ref{16}) is equal, then (\#\#) says (\ref{second}) is satisfied.
\medskip

\noindent \textbf{[Case III]} Let $\bar{c}<\bar{a}$ and $ n>1$.
We will prove that for  $i=0,1,\dots,\alpha-(b(n)+1)$
 \begin{eqnarray}
{{\alpha-i+n-1}\choose n-1}&=&d_i+\sum_{j=t_{i+1}}^{t_i-1} {{j+n-1}\choose n-1} -d_{i+1}\label{20} \\
\mbox{and} \hspace{90pt} \bar{a}&=&\bar{c}+d_{\alpha -b(n)},\label{22}
\end{eqnarray}\vspace{-2pt}
where $0\leq d_i<{{t_i+n-1}\choose n-1}$ and $t_{i+1}<t_i\leq c(n)-i+1$
together with $d_0=\bar b$, $t_0=c(n)+1$ and $t_{\alpha-b(n)}=a(n)+1$.

For $i=0$,
since ${\alpha +n-1 \choose n-1}> \bar b $, by the same way of [Case II] we have
$${{\alpha+n-1}\choose n-1}=\bar{b}+\sum_{j=t_1}^{c(n)} {{j+n-1}\choose n-1} -d_1.$$
Also, if we have equality (\ref{20}) for $i=0,1,\dots,s-1$,
then  we have ${{\alpha-s+n-1}\choose n-1} \geq  {{c(n)-(s-1)+n-1}\choose n-1}\geq {{t_{s}+n-1}\choose n-1} >d_s$.
Thus we have $t_{s+1}<t_{s}$ and we have equality (\ref{20}) for $i=s$ by the same way.
Finally, since $\bar a -\bar c < \bar a\leq { a(n) +n \choose n-1}$ by definition of $\bar a$,
we have $\bar a = \bar c +d_{\alpha -b(n)}$ and $t_{\alpha -b(n)}=a(n)+1$.
Equalities (\ref{20}) satisfies the conditions of (\#\#).
Thus by the assumption of induction of $n$, we have \begin{eqnarray}
{{\alpha-i+n-1}\choose n-1}\up{n-1} \leq d_i\up{n-1}+\sum_{j=t_{i+1}}^{t_i-1} {{j+n-1}\choose n-1}\up{n-1} -d_{i+1}\up{n-1}.\label{Y}
\end{eqnarray}
Furthermore, since  $\bar c>0$ and $d_{\alpha -b(n)}>0$ we have \begin{eqnarray}
\bar{a}\up{n-1}<\bar{c}\up{n-1}+d_{s+1}\up{n-1}.\label{Y'}
\end{eqnarray}

Then, by summating (\ref{20}) and (\ref{22}), we have the equality (\ref{15}).
By summating inequalities (\ref{Y}) and (\ref{Y'}),
we have
$$ \bigg\{ \sum_{i=b(n)+1}^\alpha {{i+n-1}\choose n-1} \bigg\}^{[+1]}+\bar{a}\up{n-1} < \bar{b}\up{n-1}+\bar{c}\up{n-1}+\bigg\{\sum_{i=a(n)+1}^{c(n)} {{i+n-1}\choose n-1}\bigg\}^{[+1]}.$$
In this case (\ref{16}) is not equal. Thus we need not consider the equality (\ref{second}).
\end{proof}
\begin{lemma}\label{rankdown}
Let $h$ and $n$ be positive integers.
Then, one has
$$h\up{n}<h\up{n+1}.$$
\end{lemma}
\begin{proof}
Let $h={{h(n+1)+n+1}\choose n+1}+\bar{h}^{(n+1)}$.
Then  $h\up{n+1}={{h(n+1)+n+1}\choose n+1}\up{n+1}+{\bar{h}^{(n+1)}}\up{n}$.
By Lemma \ref{apply}, we have
$${{h(n+1)+n+1}\choose n+1}={{n+1}\choose n+1}+\sum_{i=1}^{h(n+1)} {{i+n}\choose n}.$$
Furthermore, we have ${{n+1}\choose n+1}^{[+1]}>{{n}\choose n}^{[+1]}$.
By Lemma \ref{plus}, we have\begin{eqnarray*}
{{h(n+1)}\choose n+1}\up{n+1}+{\bar{h}^{(n+1)}}\up{n}&>&{{n}\choose n}\up{n}+\sum_{i=1}^{h(n+1)} {{i+n}\choose n}\up{n}+{\bar{h}^{(n+1)}}\up{n}\\
&\geq& \{ {{n}\choose n}+\sum_{i=1}^{h(n+1)} {{i+n}\choose n}+{\bar{h}^{(n+1)}} \}\up{n}=h\up{n},
\end{eqnarray*}
as desired.\end{proof}

\section{A combinatorial proof of persistence for monomial ideals}

Let $V$ be a set of monomials of degree $d$ and  $u=\gcd(V)$.
If $|V|>1$, we define $\main=\{v\in M^d|\ x_iu\ \mathrm{divides} \ v\} $  and $\sub =V \setminus \main$  for $i=1,2,\dots,n$.
If $ |V|=1$, then  we define \main $=V$ and $\sub=\emptyset$.
Note that if $|V|>1$, then $\sub\ne\emptyset$ and $\main \ne \emptyset$.

Before giving a combinatorial proof of persistence theorem for monomial ideals, we prepare some lemmas.
\medskip

\begin{lemma}\label{msbasic}
Let $V$ be a set of monomials of degree $d$ and $u=\gcd(V)$.
For any $i=1,2,\dots,n$, we have
\begin{eqnarray}
& &\overline{M_i}\sub \subset MV\setminus x_iV.\label{2} \\
& &|MV| \geq |\main |\up{n-1}+|\sub | \up{n-2}.  \label{(1)} 
\end{eqnarray}
Moreover, in (\ref{(1)}), the equality holds if and only if \main \  is a Gotzmann set of $K[x_1,x_2,\dots,x_n]$,
 $\frac{1}{u} \sub$  is a Gotzmann set of $K[x_1,\dots,x_{i-1},x_{i+1},\dots,x_n]$ and $x_i\sub \subset \overline{M_i}\main$.

\end{lemma}
\begin{proof}
Any element of $\overline{M_i}\sub $ can not be divided by $ux_i$.
On the other hand, $\overline{M_i}\sub \subset MV$.
Thus we have $\overline{M_i}\sub \subset MV\setminus x_iV.$

Now we have
\[
| MV|=
| M\main| +|M\sub|-|\{M\main\cap M\sub\}|.
\]
Now we have
${M\main\cap M\sub} = M\main\cap x_i\sub \subset x_i\sub$
and
$| M\sub |=| x_i\sub| +|\overline{M_i} \sub|$.
On the other hand, the inequality (\ref{mb}) says $|\overline{M_i} \sub |=|\sub|\up {n-2}$
since $\frac{1}{u} \sub \subset K[x_1,\dots,x_{i-1},x_{i+1},\dots,x_n]$.
Thus we have \begin{eqnarray*}
| MV| &\geq& |M\main| +| \overline{M_i}\sub| \\
&\geq& |\main|\up{n-1}+|\sub|\up{n-2}.
\end{eqnarray*}
Especially, equality holds if and only if \main \ and $\frac{1}{u} \sub$ are Gotzmann sets 
and $M\main\cap x_i\sub = x_i\sub$.
\end{proof}
Next  we  determine the range of $|\sub|$,
when $V$ is a Gotzmann set.

\begin{lemma}\label{hanni}
Let $V$ be a Gotzmann set of monomials of degree d.
Then, for any $i =1,2,\dots,n$, we have \begin{eqnarray}
\overline{|V|}^{(n-1)}\leq |\sub| \leq |V| \doum{n-1}.\label{A-1}
\end{eqnarray}
\end{lemma}
\begin{proof}
If $|V|=0$ or $|V|=1$, then $\overline{|V|}^{(n-1)}=|\sub|=0.$
Thus we may assume $n>1$ and $|V|>1$.
First, we consider the second inequality of (\ref{A-1}).
By Lemma  \ref{msbasic} and by the inequality (\ref{mb}), we have
$$|\sub|\up{n-2}\leq |\overline{M_i}\sub| \leq |(MV\setminus x_iV)|.$$
On the other hand, by  Lemma \ref{change}, we have
\begin{eqnarray*}|(MV\setminus x_iV)|=|MV|-|V| &=& |V|\up{n-1}-|V| \\
&=&|V|\dnm{n-1}.\end{eqnarray*}
Thus $|\sub|\up{n-2}\leq |V|\dnm{n-1}$.
Hence we have $|\sub|\leq |V|\doum{n-1}.$

We consider the first inequality of (\ref{A-1}).
If $n=2$, then $\overline{|V|}^{(n-1)}=|\sub|=1$.
Thus we may assume $n\geq3$.
Let $|V|={{a+n-1}\choose n-1}+\overline{|V|}^{(n-1)}$.
If $|\sub|<\overline{|V|}^{(n-1)}$, then
$|\main|=|V|-|\sub|>{{a+n-1}\choose n-1}.$
Thus we can write $|\main|={{a+n-1}\choose n-1}+b$ with $ b>0$.
By Lemma \ref{msbasic},  we have\begin{eqnarray*}
|MV|&\geq& |\main|\up{n-1}+|\sub|\up{n-2}\\
&=&{{a+n-1}\choose n-1}\up{n-1}+b\up{n-2}+|\sub|\up{n-2}.\end{eqnarray*}
On the other hand, by Lemma \ref{plus}, we have
\[b\up{n-2}+|\sub|\up{n-2}>\{ b+|\sub| \}\up{n-2}=(\overline{|V|}^{(n-1)})\up{n-2}.\]
Thus we have
$$|MV|>{{a+n-1}\choose n-1}\up{n-1}+(\overline{|V|}^{(n-1)})\up{n-2}=|V|\up{n-1}.$$
This is a contradiction since $V$ is a Gotzmann set.
\end{proof}

Now, we finished all preparations for following lemma which proves the Persistence Theorem immediately.

\begin{lemma}\label{divide}
Let $V$ be a Gotzmann set of monomials  of degree $d$ with $\gcd(V)=1$ and $V\ne M^d$.
Then there exist $i\in\{1,2,\dots,n\}$ which satisfies followings:
\begin{itemize}
\item[(i)] \main \  is a Gotzmann set of $K[x_1,\dots,x_n]$,
\sub \ is a Gotzmann set of $K[x_1,\dots,x_{i-1},x_{i+1},\dots,x_n]$  and $|\sub |<|V|\doum {n-1}$;\smallskip
\item[(ii)]   $x_i\sub \subset \overline{M_i}\main$;\smallskip
\item[(iii)] $\{|\main| \up{n-1} \}\up{n-1}+\{|\sub| \up{n-2}\}\up{n-2}=\{| V| \up{n-1}\}\up{n-1}$.
\end{itemize}
\end{lemma}
\begin{proof}
Now, we set $|V|=a=\sum_{j=p}^{n-1} {{a(j)+j}\choose j}$,
$|\sub| =b=\sum_{j=q}^{n-2} {{b(j)+j}\choose j}$
and $|\main|=c=\sum_{j=r}^{n-1} {{c(j)+j}\choose j}$ be the binomial representations.
Set $V\ne \emptyset$.
\medskip

\noindent \textbf{[Case(A)]} Let $|V|=1$ or $n=1$.

If $|V|=1$, then $V=M^0$ since $\gcd(V)=1$.
If $n=1$, then  $|V|=1$.
Thus we may assume $|V|>1$ and $n>1$.
\medskip

By Lemma \ref{msbasic}, if $a\up{n-1}\leq b\up{n-2}+c\up{n-1}$,
then $a\up{n-1}= b\up{n-2}+c\up{n-1}$ and $\main$ and $\sub$ are Gotzmann sets.
Thus conditions  (i) and (ii)  are satisfied.
In  {[Case(B)]} and {[Case(C)]} , we will prove that if $b<a\doum{n-1}$ then  $a\up{n-1}\leq b\up{n-2}+c\up{n-1}$ .

If $b<a\doum{n-1}$, then, by Lemma \ref{b-unique},
there exists a maximal integer $t$,
such that $n-1\geq t\geq p$ and $$\ 0\leq b-\sum_{j=t+1}^{n-1}{{a(j)+j-1}\choose j-1}<{{a(t)+t-1}\choose t-1}.$$
\newpage
Let \vspace{-12pt}
\begin{eqnarray} 
a&=&\sum_{j=t+1}^{n-1} { a(j)+j \choose j} + { a(t) +t \choose t} + a',\label{ei}\\
b& =&\sum_{j=t+1}^{n-1}{{a(j)+j-1}\choose j-1}+b',\label{bii}\\
\mbox{and} \quad c&=&a-b=\sum_{j=t+1}^{n-1}{{a(j)+j-1}\choose j}+c'.\label{sii}
\end{eqnarray}
Since $0\leq b'<{{a(t)+t-1}\choose t-1}$, we have ${{a(t)+t-1}\choose t}<c'<{{a(t)+t+1}\choose t}$.
Also, we have\begin{eqnarray}
a\up{n-1}&=&\bigg\{\sum_{j=t+1}^{n-1}{{a(j)+j}\choose j}\bigg\}^{[+1]}+{{a(t)+t}\choose t}^{[+1]}+a'\up{t-1}\label{3} \\
\mbox{and} \quad b\up{n-2}&=&\bigg\{ \sum_{j=t+1}^{n-1} {{a(j)+j-1}\choose j-1} \bigg\}^{[+1]}+b'\up{t-1}.\label{4}\end{eqnarray}
\medskip

\noindent \textbf{[Case(B)]} Let $b<a\doum{n-1}$ and $c'<{{a(t)+t}\choose t}$.

Let $c''=c'-{{a(t)+t-1}\choose t}$.
If $b'=0$, then $c'\geq {{a(t)+t}\choose t}$.
Thus $b'>0$.
On the other hand, we have $c''>0$ since $c'>{{a(t)+t-1}\choose t}$.
Since $c''< {{a(t)+t-1}\choose t-1} $, 
$c=\sum_{j=t}^{n-1}{{a(j)+j-1}\choose j}+\{\ (t-1)\mathrm{th \ binomial \ representation \ of} \ c''\}$
is $(n-1)$th binomial representation of $c$.
Thus
$$c\up{n-1}=\bigg\{ \sum_{j=t}^{n-1}{{a(j)+j-1}\choose j} \bigg\}^{[+1]}+c''\up{t-1}.$$
Thus, by (\ref{4}), we have\begin{eqnarray}
b\up{n-2}+c\up{n-1}=\bigg\{\sum_{j=t+1}^{n-1}{{a(j)+j}\choose j}\bigg\}^{[+1]}&+&{{a(t)+t-1}\choose t}^{[+1]} \nonumber \\
& & { } +b'\up{t-1}+c''\up{t-1}. \label{5} \end{eqnarray} 
Since ${{a(t)+t}\choose t}={{a(t)+t-1}\choose t}+{{a(t)+t-1}\choose t-1}$ and
 $a=b+c$ together with (\ref{ei}), (\ref{bii}) and (\ref{sii}) say $b'+c''=a'+{{a(t)+t-1}\choose t-1}$.
Hence by Lemma \ref{plus}, Lemma \ref{combinatrix} together with $b'>0$ and $ c''>0,$ we have \begin{eqnarray}
b'\up{t-1}+c''\up{t-1}\geq a'\up{t-1}+{{a(t)+t-1}\choose t-1}\up{t-1}.\label{6} \end{eqnarray}
Thus by (\ref{3}) and (\ref{5}), we have  $a\up{n-1}\leq b\up{n-2}+c\up{n-1}$.
Furthermore, if (\ref{6}) is equal,
then  Lemma \ref{combinatrix}  says 
\[
\{b'\up{t-1}\}\up{t-1}+\{c''\up{t-1}\}\up{t-1}= \{a'\up{t-1}\}\up{t-1}+\bigg\{{{a(t)+t-1}\choose t-1}\up{t-1}\bigg\}\up{t-1}.\]
Thus we have $\{a\up{n-1}\}\up{n-1} = \{b\up{n-2}\}\up{n-2}+\{c\up{n-1}\}\up{n-1}$.
\medskip

\noindent \textbf{[Case(C)]} Let $b<a\doum{n-1}$ and $c'\geq {{a(t)+t}\choose t} $.

Let $c''=c'-{{a(t)+t}\choose t}$ and  $\alpha = $max$\{ i|a(i)=a(t)\}$.
Since $\sum_{j=t+1}^{\alpha} {{a(j)+j-1}\choose j} + {{a(t)+t} \choose t}={{a(\alpha)+\alpha} \choose \alpha}$
and $c''< {{a(t)+t}\choose t-1}\leq {{a(\alpha) +\alpha}\choose \alpha-1}$,
we have
\begin{eqnarray*}
c\up{n-1}&=&\bigg\{ \sum_{j=\alpha+1}^{n-1} {{a(j)+j-1}\choose j} \bigg\}^{[+1]} +{{a(\alpha)+\alpha}\choose \alpha}^{[+1]}+c''\up{\alpha-1}\\
&=&\bigg\{\sum_{j=t+1}^{n-1} {{a(j)+j-1}\choose j}+{{a(t)+t}\choose t} \bigg\}^{[+1]}+c''\up{\alpha-1}.
\end{eqnarray*}
Thus, by (\ref{4}), we have\begin{eqnarray}
b\up{n-2}+c\up{n-1}\!=\!\bigg\{\sum_{j=t}^{n-1}{{a(j)+j}\choose j}\bigg\}^{[+1]}+b'\up{t-1}+c''\up{\alpha-1}.\label{7}\end{eqnarray}
Since $a=b+c$ together with (\ref{ei}), (\ref{bii}) and (\ref{sii}), we have $a'=b'+c''$.
By Lemmas \ref{plus} and \ref{rankdown}, we have\begin{eqnarray}
a'\up{t-1}\leq b'\up{t-1}+c''\up{t-1}\leq b'\up{t-1}+c''\up{\alpha-1}.\label{8}\end{eqnarray}
Hence by (\ref{3}) and (\ref{7}) we have $a\up{n-1}\leq b\up{n-2}+c\up{n-1}$.
Furthermore, if the inequality (\ref{8}) is equal,
then  $c'=0$ or $b'=0$ and $\alpha=t$.
In each case, we have $\{ a' \up{t-1} \} \up{t-1} = \{ {b'} \up{t-1} \} \up{t-1}+\{{c''}\up{\alpha -1}\}\up{\alpha-1}$.
Hence we have $\{a\up{n-1}\}\up{n-1} = \{b\up{n-2}\}\up{n-2}+\{c\up{n-1}\}\up{n-1}$.
\medskip

\noindent \textbf{[Case(D)]} Let $b=a\doum{n-1}$.

By Lemma \ref{msbasic}, we have $\overline{M_i}\sub\subset MV\setminus x_iV$.
But, by  (\ref{mb}) we have $|\overline{M_i}\sub| \geq b\up{n-2}$.
Now, we have $a\dnm{n-1}=a\up{n-1}-a=|(MV\setminus x_iV)|$ and $b\up{n-2}=a\dnm{n-1}$.
Thus we have $\overline{M_i}\sub = MV\setminus x_iV$.
\medskip

By [Case(B)] and [Case(C)],
if $|\sub |< a \doum {n-1}$ for some $i$,
then we have conditions (i), (ii) and (iii).
Finally, we will prove that if $|\sub|=a\doum{n-1}$ for  $i=1,2,\dots,n$ , then $V=M^d$ or $V=\emptyset$.
In [Case(D)], we see $\overline{M_i}\sub = MV\setminus x_iV$ if $|\sub|=a \doum{n-1}$.
We  claim (\#).
\begin{itemize}
\item[\textbf{(\#)}] 
Assume $|\sub|=a \doum{n-1}$ for $i=1,2,\dots,n$.
If there exist a monomial $v \in M^d$ such that $v\notin V$,
then for any $x_j$ and $x_i$ with $x_i | v$, one has $\frac{x_j}{x_i}v\notin V$.
\end{itemize}
To see why (\#) is true, we assume that $v\notin V$ and there exist $x_i$ and $x_j$ such that $\frac{x_j}{x_i}v\in V$.
Since $v\notin V$, we have $x_jv\notin x_jV$.
Thus we have $x_i \frac{x_j}{x_i}v=x_jv\in MV\setminus x_jV=\overline{M_j}\subj $.
But any element in $\overline{M_j}\subj $ does not contain $x_j$ since $\gcd(V)=1$,
this is a contradiction.
\medskip

By using (\#), if there exists a monomial $v\in M^d$ such that $v\notin V$,
then all monomials in $M^d$ do not belong to $V$.
Hence we have $V=M^d$  or $ V=\emptyset$ .
\end{proof}

We are now in the position to finish our combinatorial proof of  persistence theorem for monomial ideals.
\begin{proof}[Proof of persistence theorem for monomial ideals]
What we have to prove is  that
if $V$ is a Gotzmann set then 
 $MV$ is also a Gotzmann set.
 
Let $V$ be a Gotzmann set of degree $d$.
We use induction on $|V|$.
Firstly, for any monomial $u\in R$,
$V$ is a Gotzmann set if and only if $uV$ is a Gotzmann set
since $|V|=|uV|$ and $|MV|=|uMV|$.
Thus we may assume $\gcd(V)=1$.

If $V=M^d$, then $MV$ is also a Gotzmann set.
If $|V|=1$ then $V=M^0$.

If $V\ne M^d$ and $|V|>1$.
Lemma \ref{divide} (ii) says there exists $i\in\{1,2,\dots,n\}$ such that
$\overline{M_i}\main \supset x_i \sub$ and
$M^2\main\supset \overline{M_i}^2\main\supset x_i \overline{M_i}  \sub$.
Thus $| MV|=| M\main| +| \overline{M_i}\sub|$ and
$ |M^2V|=|M^2\main |+ | \overline{M_i}^2 \sub|.$
By  Lemma \ref{divide} (i)  and by assumption of induction,
both $M\main$ and $\overline{M_i}\sub$ are Gotzmann sets.
Hence by Lemma \ref{divide} (iii),
we have \begin{eqnarray*}
| M^2V|&=&|M^2\main| + |\overline{M_i}^2 \sub| \\
&=& \{|\main|\up {n-1}\} \up{n-1} + \{ |\sub| \up {n-2} \}\up {n-2} \\
&=&\{ |V|\up{n-1}\}\up{n-1} \\
&=& \{|MV|\}\up{n-1}.
\end{eqnarray*}
This completes the proof.
\end{proof}

\bigskip

\noindent
Satoshi Murai\\
Department of Pure and Applied Mathematics\\
Graduate School of Information Science and Technology\\
Osaka University\\
Toyonaka, Osaka, 560-0043, Japan\\
E-mail:s-murai@cr.math.sci.osaka-u.ac.jp

\end{document}